\documentclass[12pt]{article}
\usepackage{amsfonts}
\usepackage{amssymb}
\usepackage{eucal}
\usepackage{amsmath}

\begin{document}
\catcode`ð=\active
\defð{\u{g}}
\catcode`Ð=\active
\defÐ{\u{G}}
\catcode`Ý=\active
\defÝ{\. I}
\catcode`ö=\active
\defö{\"{o}}
\catcode`Ö=\active
\defÖ{\"O}
\catcode`ü=\active
\defü{\"{u}}
\catcode`Ü=\active
\defÜ{\"{U}}
\catcode`Þ=\active
\defÞ{\c{S}}
\catcode`þ=\active
\defþ{\c{s}}
\catcode`ý=\active
\defý{{\i}}
\catcode`ç=\active
\defç{\d{c}}
\catcode`Ç=\active
\defÇ{\d{C}}
\def\oast{\mathbin{\bigcirc\llap{$\ast$\kern.1cm}}}
\def\Nlim{\mathop{\rm N\!-\!lim}}
\def\nLOP{\mathop{\rm lim}}
\def\ne{\nLOP\limits_{\varepsilon\to 0}}
\def\nn{\nLOP\limits_{n\to\infty}}
\def\NLOP{\mathop{\rm N\!-\!lim}}
\def\Nn{\NLOP\limits_{n\to\infty}}
\def\Ne{\NLOP\limits_{\varepsilon\to 0}}
\def\Nm{\NLOP\limits_{m\to\infty}}
\def\Nv{\NLOP\limits_{\nu\to\infty}}
\def\cosec{\mathop{\rm cosec}\nolimits}
\def\coshec{\mathop{\rm coshec}\nolimits} \def\Ci{\mathop{\rm Ci}\nolimits}
\def\Si{\mathop{\rm Si}\nolimits} \def\oot{\mathop{\textstyle {1\over
2}}\nolimits} \def\fr#1#2{{\textstyle {#1\over#2}}}
\def\sgn{\mathop{\rm sgn}\nolimits} \def\itr{\mathop{\rm
int}\nolimits} \def \tast{\mathbin{\bigtriangleup \!\!\!\!*\,}}
\def\bast{\mathbin{\vcenter{\hrule\hbox to 8.3pt{\vrule\kern0.9pt\vbox{\hbox
{\mathstrut$\ast$}}\kern0.9pt\vrule}\hrule}}}
\def\lamb{\mathop{\lambda\!\!\!\!\lambda}}
\newcommand{\bsym}[1]{\mbox{\boldmath$#1$}}
\def\mub{\mathop{\mu\!\!\!\!\mu}}
\def\muf{\mathop{\mu\!\!\!\!\mu}}
\newcommand{\ol}{\overline}
\newcommand{\ul}{\underline}
\newcommand{\bc}{\begin{center}}
\newcommand{\ec}{\end{center}}
\newcommand{\ts}{\textstyle}
\newcommand{\ds}{\displaystyle}
\newcommand{\be}{\begin{equation}}
\newcommand{\ee}{\end{equation}}
\newcommand{\beq}{\begin{eqnarray}}
\newcommand{\eeq}{\end{eqnarray}}
\newcommand{\beqa}{\begin{eqnarray*}}
\newcommand{\eeqa}{\end{eqnarray*}}
\newcommand{\lam}{\lambda}
\def\lamb{\mathop{\lambda\!\!\!\!\lambda}}
\newcommand{\vphi}{\varphi}
\newcommand{\Aa}{{\cal A}}
\newcommand{\BB}{{\cal B}}
\newcommand{\DD}{{\cal D}}
\newcommand{\EE}{{\cal E}}
\newcommand{\GG}{{\cal G}}
\newcommand{\RR}{{\cal R}}
\newcommand{\TT}{{\cal T}}
\newcommand{\XX}{{\cal X}}
\newcommand{\YY}{{\cal Y}}
\newcommand{\ZZ}{{\cal Z}}
\newcommand{\CC}{{\cal C}}
\newcommand{\n}{\nu}
\newcommand{\la}{\langle}
\newcommand{\ra}{\rangle}
\newcommand{\bcup}{\textstyle{\bigcup}}
\newcommand{\bcap}{\textstyle{\bigcap}}
\newcommand{\reals}{{\mathbb R}}
\newcommand{\integers}{{\mathbb Z}}
\newcommand{\naturals}{{\mathbb N}}

\medskip
\begin{center}\LARGE {\bf Defining Incomplete Gamma Type Function with Negative
Arguments and Polygamma functions $\psi^{(n)}(-m)$}
\end{center}

\medskip
\begin{center}{\bf Emin Özçað and Ýnci Ege
}
\end{center}

\bigskip \noindent
{\bf Abstract.} {\footnotesize In this paper the incomplete gamma
function $\gamma(\alpha,x)$ and its derivative is considered for
negative values of $\alpha $ and the incomplete gamma type
function $\gamma_*(\alpha,x_-)$ is introduced. Further the
polygamma functions $\psi^{(n)}(x)$ are defined for negative
integers via the neutrix setting.}

\bigskip \noindent {\bf {\underline {AMS Mathematics Subject
Classification (2010)}}:} 33B15, 33B20, 33D05

\bigskip \noindent {\bf {\underline {Key words and phrases}}:}
Gamma function, incomplete gamma function, Digamma Function,
Polygamma Function, neutrix, neutrix limit.

\bigskip \noindent
\bc {\Large {\bf 1. Introduction}} \ec

\bigskip
\noindent The incomplete gamma function of real variable and its
complement are defined via integrals \beq \gamma(\alpha,x) =
\int_0^x u^{\alpha - 1}e^{-u} \,du \eeq and \beq \Gamma(\alpha,x)
= \int_x^\infty u^{\alpha - 1}e^{-u} \,du \eeq respectively, for
$\alpha > 0$ and $x\geq 0$ see \cite{kn:Abr,kn:rr}.

\noindent Note that the definition of $\gamma(\alpha,x)$ does not
valid  for negative $\alpha ,$ but it can be extended by the
identity $$ \gamma(\alpha,x) = \Gamma(\alpha )-\Gamma(\alpha,x).$$

\noindent The incomplete gamma function does not exists for
negative  integers $\alpha $ or zero \cite{kn:Thomp1}. However,
Fisher et. al. \cite{kn:fbkil} defined $\gamma (0,x)$ by
$$\gamma(0,x) = \int_0^x u^{-1} (e^{-u}-1) \,du+\ln x .$$
The problem of evaluating incomplete gamma function is subject of
some earlier articles. The general problem in which $\alpha $ and
$x$ are complex was considered by Winitzki in \cite{kn:Win}, but
no method is found to be satisfactory in cases where $Re [x]<0.$
The incomplete gamma function with negative arguments are
difficult to compute, see \cite{kn:gautschi}. In \cite{kn:Thomp2}
Thompson gave the algorithm for accurately computing the
incomplete gamma function $\gamma (\alpha, x)$ in the cases where
$\alpha =n+1/2, n\in \integers $ and $x<0.$

\bigskip \noindent
\bc {\Large {\bf 2. Incomplete Gamma Function}} \ec

\bigskip \noindent
A series expansion of incomplete gamma function can be obtained by
replacing the exponential in (1) with its Maclaurin series, the
result is that \be \gamma(\alpha,x)=\sum_{k=0}^\infty
\frac{(-1)^k}{k!(\alpha +k)} x^{\alpha +k}. \ee A single
integration by parts in (1) yields the recurrence relation \be
\gamma (\alpha+1,x) = \alpha \gamma(\alpha,x) -x^{\alpha}e^{-x}
\ee and thus the incomplete gamma function can be defined for
$\alpha<0$ and $\alpha \neq -1,-2,\ldots.$

\medskip \noindent
By repeatedly applying recurrence relation, we obtain $$
\gamma(\alpha,x)= x^{\alpha -1}e^{-x}\sum_{k=0}^\infty
\frac{x^k}{(\alpha)_k}+\frac{1}{(\alpha)_m}\gamma(\alpha +m ,x) $$
where $(.)_k$ is Pochhamer's symbol \cite{kn:Abr}. The last term
on the right-hand side disappears in the limit $m\rightarrow
\infty.$

\medskip \noindent
By regularization  we have \beq \gamma(\alpha,x) = \int_0^x
u^{\alpha -1} \Bigl [ e^{-u} -\sum_{i=0}^{m-1} \frac{(-u)^i}{i!}
\Bigr ]\,du +\sum_{i=0}^{m-1}\frac{(-1)^i x^{\alpha +i}} {(\alpha
+ i)i!}. \eeq for $-m < \alpha < -m + 1 $ and $x>0.$ It follows
from the definition of gamma function that
$$\lim_{x\to\infty}\gamma(\alpha,x) = \Gamma (\alpha) $$ for
$\alpha \neq 0,-1,-2, \ldots,$ see \cite{kn:gs,kn:ozinci}.

\medskip \noindent In the following we let N be the neutrix
\cite{kn:c,kn:f1,kn:fbkil} having domain $N'
={\{\varepsilon:\,0<\varepsilon <\infty}\}$ and range $N''$ the
real numbers, with negligible functions finite linear sums of the
functions \be \varepsilon^{\lambda} \ln^{r-1} \varepsilon, \quad
\ln^r \varepsilon \qquad (\lambda <0, \quad  r \in \integers
^+)\ee and all functions of $\varepsilon$ which converge to zero
in the normal sense as $\varepsilon$ tends to zero.

\medskip \noindent
If $f(\varepsilon)$ is a real (or complex) valued function defined
on $N'$ and if it is possible to find a constant $c$ such that
$f(\varepsilon) -c$ is in N, then $c$ is called the neutrix limit
of $f(\varepsilon)$ as $\varepsilon \rightarrow 0$ and we write $
\Nlim_{\varepsilon\to 0} f(\varepsilon) =c.$

\medskip \noindent
Note that if a function $f(\varepsilon)$ tends to $c$ in the
normal sense as $\varepsilon$ tends to zero, it converges to $c$
in the neutrix sense.

\medskip \noindent  On using
equation (5), the incomplete gamma function $\gamma(\alpha,x)$ was
also defined by \be \gamma(\alpha,x) = \Ne \int_{\varepsilon}^x
u^{\alpha -1} e^{-u}\,du \ee for all $\alpha\in\reals$ and $x>0,$
and it was shown that $\lim_{x\to\infty}\gamma(-m,x) = \Gamma
(-m)$ for $m \in \naturals,$ see \cite{kn:fbkil}.

\medskip \noindent
Further, the r-th derivative of  $\gamma (\alpha,x)$ was similarly
defined by \beq \gamma^{(r)}(\alpha,x) = \Ne \int_{\varepsilon}^x
u^{\alpha -1} \ln^r u e^{-u}\, du \eeq for all $\alpha $ and
$r=0,1,2,\ldots,$ provided that the neutrix limit exists, see
\cite{kn:ozinci}.

\medskip \noindent
By using the neutrix N given in equation(6), the interesting
formula
$$ \gamma(-m,x) +\frac{1}{m} \gamma(-m+1,x) =\frac{(-1)^m}{m m!}
-\frac{1}{m} e^{-x}x^{-m}
$$
was obtained in \cite{kn:ozinci} for $m\in \naturals .$

\medskip \noindent
It can be easily seen that equation (3) is also valid for $\alpha
< 0$ and $\alpha \neq -1,-2,\ldots.$ In fact it follows from the
definition that \beqa \gamma (\alpha,x) &=&\Ne
\int_{\varepsilon}^x u^{\alpha -1}e^{-u}\,du=\Ne \Bigl [
\sum_{k=0}^\infty \frac{(-1)^k}{k!}
\int_\varepsilon ^x u^{\alpha+k-1}\,du \Bigr ] \\
&=& \Ne \sum_{k=0}^\infty \frac{(-1)^k}{k!(\alpha+k)}\bigl
[x^{\alpha+k}-\varepsilon^{\alpha+k} \bigr ] \\
&=& \sum_{k=0}^\infty
\frac{(-1)^k}{k!(\alpha+k)}x^{\alpha+k}.\eeqa

\medskip \noindent
Now assume that $\alpha\in \integers^-,$ then writing \beqa
\int_{\varepsilon}^x u^{-m-1} e^{-u}\,du &=& \sum_{\substack{k=0,
\\k\neq m}}^\infty \frac{(-1)^k}{k!}\int_{\varepsilon}^x
u^{-m+k-1}\,du+\frac{(-1)^m}{m!}\int_{\varepsilon}^x u^{-1}\,du \hskip0.9in \\
&=& \sum_{\substack{k=0, \\ k\neq m}}^\infty
\frac{(-1)^k}{k!(k-m)}x^{k-m}+\frac{(-1)^m}{m!}\ln x+
O(\varepsilon) \eeqa where $O(\varepsilon)$ is negligible
function, and taking the neutrix limit we obtain
$$\gamma(-m,x)=\sum_{\substack{k=0, \\k\neq m}}^\infty\frac
{(-1)^k}{k!(k-m)}x^{k-m}+\frac{(-1)^m}{m!}\ln x $$ for $x>0$ and
$m\in \naturals .$ If $m=0$ we particularly have
$$ \gamma (0,x)=\sum_{k=1}^\infty
\frac{(-x)^k}{kk!}+\ln x$$ and thus it follows by the calculation
of Mathematica that
$$\gamma(0,1/2)\approx -1.13699, \quad \gamma(0, 3/4)\approx -0.917556, \quad \gamma(0,
1)\approx -0.7966.$$

\medskip \noindent
Similar equation for $\gamma^{(r)}(\alpha,x)$ can be obtained
using the following lemma.

\medskip \noindent
{\bf Lemma 2.1}  \beqa \Ne \int_\varepsilon^x u^{-1}\ln ^r
u\,du &=& \frac{\ln^{r+1} x}{r+1} \\
\Ne \int_\varepsilon^x u^\alpha \ln ^r u\,du &=& \sum_{k=0}^{r-1}
\frac{(-1)^k r! x^{\alpha+1}\ln ^{r-k} x}{(\alpha+1)^{k+1}(r-k)!}+
\frac{(-1)^r r! x^{\alpha+1}}{(\alpha+1)^{r+1}} \eeqa for all
$\alpha \neq -1, x>0$ and $r\in \naturals.$

\medskip \noindent
{\bf Proof.} Straight forward.

\bigskip \noindent Now let $\alpha <0$ and $\alpha \neq
0,-1,-2,\ldots.$ Then it follows from lemma 2.1 that \beqa
\gamma^{(r)}(\alpha,x)&=& \Ne
\int_\varepsilon ^x u^{\alpha-1}\ln ^{r} u e^{-u}\,du \\
&=& \Ne \sum_{k=0}^\infty \frac{(-1)^k}{k!} \int_\varepsilon ^x
u^{\alpha+k-1} \ln^{r} u\,du  \\
&=& \sum_{k=0}^\infty \sum_{i=0}^{r-1}\frac{(-1)^{i+k}
r!}{k!(r-i)!(\alpha +k)^{i+1}}x^{\alpha +k}\ln ^{r-i}x +
\\
&&\hskip0.9in +\sum_{k=0}^\infty \frac{(-1)^{k+r} r!}{k!(\alpha
+k)^{r+1}} x^{\alpha +k} \eeqa and similarly \beqa
\gamma^{(r)}(-m,x)&=&
\Ne \int_\varepsilon ^x u^{-m-1}\ln ^{r} u e^{-u}\,du \\
&=& \sum_{\substack{k=0, \\ k\neq m}}^\infty
\sum_{i=0}^{r-1}\frac{(-1)^{i+k}
r!}{k!(r-i)!(m-k)^{i+1}}x^{k-m}\ln ^{r-i}x + \\
&&\hskip0.4in +\sum_{\substack{k=0, \\ k\neq m}}^\infty
\frac{(-1)^{k+r} r!}{k!(m-k)^{r+1}}
x^{k-m}+\frac{(-1)^m}{(r+1)m!}\ln^{r+1} x \eeqa for $r=1,2,\ldots$
and $m\in \naturals.$

\bigskip \noindent
It was indicated in \cite{kn:f1} that equation (1) could be
replaced by the equation \be \gamma(\alpha,x) = \int_0^x
|u|^{\alpha - 1}e^{-u} \,du \ee and this equation was used to
define $\gamma(\alpha,x)$ for all $x$ and $\alpha >0,$ the
integral again diverging for $\alpha\leq 0.$

\medskip \noindent
The locally summable function $\gamma_*(\alpha,x_-)$ for $\alpha
>0$ was then defined by
\beq \gamma_*(\alpha,x_-)&=& \left \{ \begin{array}{cc}
\int_0^x |u|^{\alpha-1}e^{-u}\,du, & x\leq 0, \\
0, & x>0
\end{array} \right.  \nonumber \\
&=& \int_0^{-x_-} |u|^{\alpha-1}e^{-u}\,du. \eeq If $\alpha >0,$
then the recurrence relation  \be \gamma_*(\alpha +1,x_-)=-\alpha
\gamma_*(\alpha ,x_-)-x_-^\alpha e^{-x}. \ee holds. So we can use
equation (11) to extend the definition of $\gamma_*(\alpha ,x_-)$
to negative non-integer values of $\alpha.$

\medskip \noindent
More generally it can be easily proved that \be
\gamma_*(\alpha,x_-) = \int_0^{-x_-} |u|^{\alpha -1} \Bigl [
e^{-u} -\sum_{i=0}^{m-1} \frac{(-u)^i}{i!} \Bigr ]\,du
-\sum_{i=0}^{m-1}\frac{x_-^{\alpha +i}} {(\alpha + i)i!}. \ee if
$-m<\alpha <-m+1$ for $m=1,2,\ldots.$

\medskip \noindent
Now if $-m<\alpha <-m+1, m\in \naturals $ and $x<0,$ we have
$$\int_{-\varepsilon}^{-x_-} |u|^{\alpha -1} e^{-u}\,du =
\int_{-\varepsilon}^{-x_-} |u|^{\alpha -1} \Bigl [ e^{-u}
-\sum_{i=0}^{m-1} \frac{(-u)^i}{i!} \Bigr ]\,du -\sum_{i=0}^{m-1}
\frac{\Bigl [|x|^{\alpha+i}-\varepsilon^{\alpha+i}\Bigr
]}{(\alpha+i)i!}$$ and thus
$$\gamma_*(\alpha ,x_-) =\Ne \int_{-\varepsilon}^{-x_-} |u|^{\alpha -1}
e^{-u}\,du $$ on using equation (12). This suggests that we define
$\gamma_*(-m,x_-)$ by \be \gamma_*(-m,x_-) =\Ne
\int_{-\varepsilon}^{-x_-} |u|^{-m-1} e^{-u}\,du \ee for $x<0$ and
$m\in \naturals .$

\medskip \noindent
If $x<0,$ then we simply write
$$\int_{-\varepsilon}^{-x_-} |u|^{-1} e^{-u}\,du =
\int_{-\varepsilon}^{-x_-} |u|^{-1} (e^{-u}-1)\,du-\ln |x|+\ln
\varepsilon $$ thus we have \beqa \Ne \int_{-\varepsilon}^{-x_-}
|u|^{-1} e^{-u}\,du
&=& \int_0^{-x_-} |u|^{-1} (e^{-u}-1)\,du-\ln x_- \\
&=& \gamma_*(0,x_-). \eeqa Next, writing similarly \beqa
\int_{-\varepsilon}^{-x_-} |u|^{-m-1} e^{-u}\,du &=&
\int_{-\varepsilon}^{-x_-} |u|^{-m-1} \Bigl [ e^{-u} -\sum_{i=0}^m
\frac{(-u)^i}{i!} \Bigr ]\,du\hskip0.9in
\\
&& \hskip0.2in -\sum_{i=0}^{m-1} \frac{\Bigl
[|x|^{i-m}-\varepsilon^{i-m}\Bigr ]}{(m-i)i!}-\frac{1}{m!}\Bigl
[\ln |x|-\ln \varepsilon \Bigr ] \eeqa for $x<0$ and $m\in
\naturals,$ then it follows that \beq &&\gamma_*(-m ,x_-) = \Ne
\int_{-\varepsilon}^{-x_-} |u|^{-m-1}
e^{-u}\,du \hskip2.7in \nonumber \\
&&\hskip0.4in =\int_0^{-x_-} |u|^{-m-1} \Bigl [ e^{-u}
-\sum_{i=0}^m \frac{(-u)^i}{i!} \Bigr ]\,du-\sum_{i=0}^{m-1}
\frac{x_-^{i-m}}{(m-i)i!}-\frac{1}{m!}\ln x_- .\eeq

\medskip \noindent
Equation (14) can be regarded as the definition of incomplete
gamma function $\gamma_*(-m ,x_-)$ for negative integers and also
written in the form \beq \gamma_*(-m ,x_-) &=& \int_0^{-1}
|u|^{-m-1} \Bigl [ e^{-u}
-\sum_{i=0}^m \frac{(-u)^i}{i!} \Bigr ]\,du \nonumber \\
&&\hskip0.4in  +\int_{-1}^{-x_-} |u|^{-m-1}
e^{-u}\,du+\sum_{i=0}^{m-1}\frac{1}{(m-i)i!}\eeq Replacing
$e^{-x}$ by its Maclaurin series yields in equation (13) \be
\gamma_*(-m ,x_-)=\sum_{\substack{k=0, \\ k\neq m}}^\infty
\frac{x_-^{k-m}}{(m-k)k!}-\frac{1}{m!}\ln x_-\quad (m\in
\naturals). \ee If $m=0$ we particular have \be
\gamma_*(0,x_-)=-\sum_{k=1}^\infty \frac{x_-^k}{kk!}-\ln x_-. \ee
and by Mathematica
\beqa &&\gamma_*(0 ,-1/3)\approx 0.735308 \quad
\gamma_*(0,-1/2)\approx 0.122996 \\ && \gamma_*(0, -3/4)\approx
-0.630117
 \quad \gamma_*(0,-1) \approx -1.3179.\eeqa
 Differentiating equation (12) we get \beqa
\gamma_*'(\alpha,x_-) &=& \int_0^{-x_-} |u|^{\alpha -1}\ln |u|
\Bigl [e^{-u} -\sum_{i=0}^{m-1} \frac{(-u)^i}{i!} \Bigr ]\,du \\
&&\hskip0.4in -\sum_{i=0}^{m-1}\frac{(\alpha+i)x_-^{\alpha +i}\ln
x_- -x_-^{\alpha+i}} {(\alpha + i)^2i!} \eeqa for $-m<\alpha
<-m+1$ and $x<0.$

\medskip \noindent
On the other hand, arranging the following integral and taking the
neutrix limit \beqa \int_{-\varepsilon}^{-x_-} |u|^{\alpha -1}\ln
|u| e^{-u}\,du &=& \int_{-\varepsilon}^{-x_-} |u|^{\alpha -1}\ln
|u| \Bigl [ e^{-u}
-\sum_{i=0}^{m-1}\frac{(-u)^i}{i!} \Bigr ]\,du \\
&&\hskip0.4in
+\sum_{i=0}^{m-1}\frac{1}{i!}\int_{-\varepsilon}^{-x_-}
|u|^{\alpha +i-1}\ln |u| \,du \\
&=& \int_{-\varepsilon}^{-x_-} |u|^{\alpha -1}\ln |u| \Bigl
[e^{-u}
-\sum_{i=0}^{m-1}\frac{(-u)^i}{i!} \Bigr ]\,du \\
&&\hskip0.1in -\sum_{i=0}^{m-1} \frac{1}{i!(\alpha + i)} \Bigl
[x_-^{\alpha+i}\ln x_- -\varepsilon^{\alpha+i}\ln
\varepsilon \Bigr ] \hskip0.4in \\
&&\hskip0.4in +\sum_{i=0}^{m-1}\frac{1}{i!(\alpha +i)^2} \Bigl
[x_-^{\alpha+i} -\varepsilon^{\alpha+i}\Bigr ]. \eeqa Thus
$$\gamma_*'(\alpha,x_-) =\Ne \int_{-\varepsilon}^{-x_-} |u|^{\alpha -1}\ln |u|
 e^{-u}\,du$$ for $\alpha\neq 0,-1,-2,\ldots$ and $x<0.$

\medskip \noindent
More generally it can be shown that
\be \gamma_*^{(r)}(\alpha,x_-)
=\Ne \int_{-\varepsilon}^{-x_-} |u|^{\alpha -1}\ln^r |u|
e^{-u}\,du\ee for $\alpha\neq 0,-1,-2,\ldots, r=0,1,2,\ldots$ and
$x<0.$ This suggests the following definition.

\medskip \noindent
{\bf Definition 2.2} The r-th derivative of incomplete gamma
function $\gamma_*^{(r)}(-m,x_-)$ is defined by \be
\gamma_*^{(r)}(-m,x_-) =\Ne \int_{-\varepsilon}^{-x_-} |u|^{-m
-1}\ln^r |u| e^{-u}\,du \ee for $r,m=0,1,2,\ldots$ and $x<0$
provided that the neutrix limit exists.

\medskip \noindent
Equation(18) will then define $\gamma_*^{(r)}(\alpha,x_-)$ for all
$\alpha $ and $r=0,1,2,\ldots.$ To prove the neutrix limit above
exists we need the following lemma.

\medskip \noindent
{\bf Lemma 2.3} \beq \int \ln^r |u| \,du =
\sum_{i=0}^{r-1}\frac{(-1)^ir!}{(r-i)!} u\ln^{r-i}|u| + (-1)^rr!u
\eeq and \beq \int u^{-s-1}\ln^r |u| \,du =
-\sum_{i=0}^{r-1}\frac{r!}{(r-i)!} s^{-i-1} u^{-s} \ln^{r-i} |u|
-r! s^{-r-1}u^{-s} \eeq for $r,s =1,2,\ldots.$

\medskip \noindent
{\bf Proof. } Equation (20) follows by induction and equation (21)
follows on using equation (20) and making the substitution
$w=u^{-s}.$

\bigskip \noindent
{\bf Theorem 2.4} The functions $\gamma_*^{(r)}(0,x_-)$ and
$\gamma_*^{(r)}(-m,x_-)$ exist for $x<0$ and
\be
\gamma_*^{(r)}(0,x_-) =\int_{-1}^{-x_-} |u|^{-1}\ln^r |u|
e^{-u}\,du +\int_0^{-1} |u|^{-1}\ln^r |u|[e^{-u} -1] \,du \ee for
$r= 0,1,2,\ldots$ and \beq &&\gamma_*^{(r)}(-m,x_-) =
\int_{-1}^{-x_-} |u|^{-m -1} \ln^r |u| e^{-u}\,du +\hskip2.5in
\nonumber \\
&&\hskip0.5in +\int_0^{-1} |u|^{-m -1}\ln^r |u| \Bigl [ e^{-u}
-\sum_{i=0}^m \frac{(-u)^i}{i!} \Bigr ]\,du
+\sum_{i=0}^{m-1}\frac{r!(m-i)^{-r-1}}{i!} \eeq for $r,m
=1,2,\ldots.$

\medskip\noindent
{\bf Proof. } We have \beqa
&&\hskip0.2in \int_{-\varepsilon}^{-x_-} |u|^{-1}\ln^r |u| e^{-u} \,du = \hskip3.9in \\
&=& \int_{-1}^{-x_-} |u|^{-1}\ln^r |u| e^{-u}\,du
+\int_{-\varepsilon}^{-1} |u|^{-1}\ln^r |u|
[e^{-u} -1] \,du +\int_{-\varepsilon}^{-1} |u|^{-1}\ln^r |u|\,du \\
&=& \int_{-1}^{-x_-} |u|^{-1}\ln^r |u| e^{-u} \,du
+\int_{-\varepsilon}^{-1} |u|^{-1}\ln^r |u| [e^{-u} -1]\,du
+\frac{\ln^{r+1}\varepsilon}{r+1}. \eeqa for $r=0,1,2,\ldots.$ It
follows that
\beqa
&& \Nlim_{\varepsilon\to 0}\int_{-\varepsilon}^{-x_-} |u|^{-1}\ln^r |u| e^{-u} \,du =\hskip3.0in \\
&&\hskip1.0in  =\int_{-1}^{-x_-} |u|^{-1}\ln^r |u| e^{-u} \,du
+\int_0^{-1} |u|^{-1}\ln^r |u|[e^{-u} -1] \,du. \eeqa and so
equation (22) follows.

\medskip \noindent
Now let us consider
\beqa &&
\int_{-\varepsilon}^{-x_-} |u|^{-m - 1}\ln^r |u| \,e^{-u}\,du
=\int_{-1}^{-x_-} |u|^{-m-1}\ln^r
|u| e^{-u} \,du + \hskip3.3in \\
&&\hskip2.0in +\int_{-\varepsilon}^{-1} |u|^{-m -1}\ln^r |u| \Bigl
[e^{-u}-\sum_{i=0}^m \frac{(-u)^i}{i!} \Bigr ]\,du
\\
&&\hskip2.5in +\sum_{i=0}^m
\frac{(-1)^i}{i!}\int_{-\varepsilon}^{-1} |u|^{-m +i-1}\ln^r
|u|\,du \eeqa \beqa && \hskip0.1in =\int_{-1}^{-x_-} |u|^{-m
-1}\ln^r |u| e^{-u} \,du +\int_{-\varepsilon}^{-1}
|u|^{-m -1}\ln^r |u|\Bigl [ e^{-u} -\sum_{i=0}^m\frac{(-u)^i}{i!} \Bigr ]\,du + \\
&&\hskip0.4in +\sum_{i=0}^{m-1}\sum_{j=0}^{r-1}\frac{(-1)^{i-m-1}
r!}{i! (r-j)!}
(m-i)^{-j-1}\varepsilon^{m-i}\ln^{r-j} \varepsilon+ \\
&&\hskip0.7in -\sum_{i=0}^{m-1}\frac{(-1)^{i-m-1}(m-i)^{-r-1}
r!}{i!}\Bigl [ (-1)^{m-i} -\varepsilon^{m-i})\Bigr
]+\frac{\ln^{r+1}\varepsilon}{m! (r+1)}. \eeqa
Thus
\beqa
&&\gamma_*^{(r)}(-m, x_-) = \Nlim_{\varepsilon\to
0}\int_{-\varepsilon}^{-x_-} |u|^{-m-1}\ln^r |u| e^{-u} \,du \hskip2.5in \\
&&\hskip0.4in = \int_{-1}^{-x_-} |u|^{-m - 1}\ln^r |u| e^{-u}\,du
+ \int_0^{-1} |u|^{-m -1}\ln^r |u| \Bigl [e^{-u} -\sum_{i=0}^m
\frac{(-u)^i}{i!} \Bigr ]\,du \\
&& \hskip0.8in  +\sum_{i=0}^{m-1}\frac{r!(m-i)^{-r-1}}{i!}. \eeqa
for $r,m =1,2,\ldots.$ Equation (23) follows.

\bigskip
\bc{\large {\bf 3. Polygamma Functions $\psi^{(n)}(-m)$}} \ec

\bigskip \noindent
The polygamma function is defined by \be
\psi^{(n)}(x)=\frac{d^n}{dx^n}\psi(x)=\frac{d^{n+1}}{dx^{n+1}}\ln
\Gamma(x) \qquad (x>0). \ee It may be represented as
$$\psi^{(n)}(x)= (-1)^{n+1}\int_0^\infty \frac{t^n
e^{-xt}}{1-e^{-t}}$$ which holds for $x>0,$ and \be \psi^{(n)}(x)=
-\int_0^1\frac{t^{x-1}\ln ^nt}{1-t}\,dt. \ee

\medskip \noindent
It satisfies the recurrence relation \be
\psi^{(n)}(x+1)=\psi^{(n)}(x)+\frac{(-1)^nn!}{x^{n+1}} \ee see
\cite{kn:rr}. This is used to define the polygamma function for
negative non-integer values of $x.$ Thus if $-m<x<-m+1, \quad
m=1,2,\ldots,$ then \be \psi^{(n)}(x)=-\int_0^1
\frac{t^{x+m-1}}{1-t}\ln^n
t\,dt-\sum_{k=0}^{m-1}\frac{(-1)^nn!}{(x+k)^{n+1}}.\ee

\medskip \noindent
Kölbig gave the formulae for the integral $\int_0^1 t^{\lambda
-1}(1-t)^{-\nu}\ln ^mt\,dt$ for integer and half-integer values of
$\lambda $ and $\nu $ in \cite{kn:kol3}. As the integral
representation of the polygamma function is similar to the
integral mentioned above, we prove the existence of the integral
in equation (25) for all values of $x$ by using the neutrix limit.

\medskip \noindent
Now we let $N$ be a neutrix having domain the open interval $\{
\epsilon : 0<\varepsilon < {1\over 2}\}$ with the same negligible
functions as in equation (6). We first of all need the following
lemma.

\medskip \noindent
{\bf Lemma 3.1.} The neutrix limits as $\varepsilon $ tends to
zero of the functions $$ \int_{\varepsilon}^{1/2} t^x \ln ^nt\ln
^r(1-t)\,dt, \qquad \int_{1/2}^{1-\varepsilon}(1-t)^x\ln ^nt\ln
^r(1-t)\,dt $$ exists for $n,r=0,1,2,\ldots$ and all $x.$

\medskip \noindent
{\bf Proof.} Suppose first of all that $n=r=0.$ Then
$$ \int_\varepsilon ^{1/2} t^x \,dt = \left \{ \begin{array}{cc}
\frac{2^{-x-1}-\varepsilon^{x+1}}{x+1}, & x\neq -1, \\
-\ln 2-\ln \varepsilon, & x=-1
\end{array} \right. $$
and so $\Ne \int_\varepsilon^{1/2} t^x \,dt $ exists for all $x.$

\medskip \noindent
Now suppose that $r=0$ and that $ \Ne
\int_\varepsilon^{1/2} t^x\ln ^n t \,dt $ exists for some
nonnegative integer n and all $x.$ Then
$$ \int_\varepsilon ^{1/2} t^x\ln^{n+1} t \,dt = \left \{ \begin{array}{cc}
\frac{-2^{-x-1}\ln^{n+1} 2-\varepsilon ^{x+1}\ln^{n+1}
\varepsilon}{x+1}- \frac{n+1}{x+1}\int_\varepsilon ^{1/2} t^x\ln^n
t \,dt, & x\neq -1,
\\ \\ \frac{(-1)^n\ln^{n+2} 2-\ln^{n+2} \varepsilon }{m+2}, & x=-1
\end{array} \right.$$
and it follows by induction that $ \Ne \int_\epsilon^{1/2} t^x
\ln^n t \,dt $ exists for $n=0,1,2,\ldots$ and all $x.$

\medskip \noindent
Finally we note that we can write
$$ \ln^r (1-t) = \sum_{i=1}^\infty \alpha_{ir} t^i $$
for $r=1,2,\ldots,$ the expansion being valid for $|t|<1.$
Choosing a positive integer k such that $x+k>-1,$ we have \beqa
&&\int_\varepsilon ^{1/2} t^x \ln^n t \ln^r(1-t) \,dt =\hskip2.8in \\
&&\hskip0.2in = \sum_{i=1}^{k-1}\alpha_{in}\int_\varepsilon ^{1/2}
t^{x+i} \ln^n t\,dt +\sum_{i=k}^\infty
\alpha_{in}\int_\varepsilon^{1/2} t^{x+i} \ln^n t\,dt.
\eeqa
It
follows from what we have just proved that
$$\Ne \sum_{i=1}^{k-1}\alpha_{in}\int_\varepsilon ^{1/2} t^{x+i}
\ln^n t\,dt $$ exists and further \beqa \Ne \sum_{i=k}^\infty
\alpha_{in}\int_\varepsilon ^{1/2} t^{x+i} \ln^n t\,dt &=& \ne
\sum_{i=k}^\infty \alpha_{in}
\int_\varepsilon ^{1/2} t^{x+i} \ln^n t\,dt \\
&=& \sum_{i=k}^\infty \alpha_{in} \int_0^{1/2} t^{x+i} \ln^n
t\,dt, \eeqa proving that the neutrix limit of $\int_\varepsilon
^{1/2} t^x \ln^n t \ln^r(1-t) \,dt$ exists for $n,r=0,1,2,\ldots$
and all $x.$ Making the substitution $1-t=u$ in
$$\int_{1/2}^{1-\varepsilon }
(1-t)^x \ln ^n t \ln ^r t\,dt,$$ it follows that
$\int_{1/2}^{1-\varepsilon } (1-t)^x \ln ^n t \ln ^r t\,dt$ also
exits for $n,r=0,1,2,\ldots$ and all $x.$

\medskip \noindent
{\bf Lemma 3.2.} The neutrix limit as $\varepsilon \rightarrow 0$
of the integral $\int_{\varepsilon }^1 t^{-m-1} \ln ^n t\,dt$
exists for $m,n=1,2,\ldots$ and \be \Ne \int_{\varepsilon }^1
t^{-m-1} \ln ^n t\,dt =-\frac{n!}{m^{n+1}}. \ee

\medskip \noindent
{\bf Proof.} Integrating by parts, we have
$$\int_{\varepsilon}^1 t^{-m-1} \ln t\,dt = m^{-1}\varepsilon^{-m} \ln \varepsilon
+m^{-1}\int_{\varepsilon}^1 t^{-m-1} \,dt $$ and so $$\Ne
\int_{\varepsilon}^1 t^{-m-1} \ln t\,dt = -\frac{1}{n^2} $$
proving equation (28) for $n=1$ and $m=1,2,\ldots.$

\medskip \noindent
Now assume that equation (28) holds for some $m$ and
$n=1,2,\ldots.$ Then \beqa \int_{\varepsilon}^1 t^{-m-2} \ln ^n
t\,dt &=& (m+1)^{-1}\varepsilon^{-m-1}\ln ^n\varepsilon
+\frac{n}{m+1}\int_{\varepsilon}^1 t^{-m-2} \ln ^{n-1} t\,dt \\
&=& (m+1)^{-1}\varepsilon^{-m-1}\ln ^n\varepsilon +
\frac{n}{m+1}\frac{-(n-1)!}{(m+1)^n} \eeqa and it follows that
$$\Ne \int_{\varepsilon}^1 t^{-m-2} \ln ^n t\,dt
=-\frac{n!}{(m+1)^{n+1}}$$ proving equation (28) for $m+1$ and
$n=1,2,\ldots.$

\medskip\noindent
Using the regularization and the neutrix limit, we prove the
following theorem.

\medskip \noindent
{\bf Theorem 3.3.} The function $\psi^{(n)}(x)$ exists for
$n=0,1,2,\ldots,$ and all $x.$

\medskip \noindent
{\bf Proof.} Choose positive integer r such that $x>-r.$ Then we
can write \beqa \int_\varepsilon ^{1-\varepsilon
}\frac{t^{x-1}}{1-t}\ln^n t\,dt &=& \int_\varepsilon ^{1/2}
t^{x-1}\ln^n t \Bigl [\frac{1}{1-t}-\sum_{i=0}^{r-1}(-1)^i t^i
\Bigr ]\,dt \\
&& \hskip0.1in + \sum_{i=0}^{r-1}(-1)^i \int_\varepsilon ^{1/2}
t^{x+i-1}\ln^n t\,dt+\int_{1/2}^{1-\varepsilon }
\frac{t^{x-1}}{1-t}\ln^n t\,dt. \eeqa

We have \beqa && \ne \int_\varepsilon
^{1/2}\frac{t^{x-1}}{1-t}\ln^n t\Bigl
[\frac{1}{1-t}-\sum_{i=0}^{r-1}(-1)^i t^i \Bigr ]\,dt =
\hskip2.0in
\\ &&\hskip1.7in =\int_0^{1/2} \frac{t^{x-1}}{1-t}\ln^n t \Bigl
[\frac{1}{1-t}-\sum_{i=0}^{r-1}(-1)^i t^i \Bigr ]\,dt \eeqa and
$$\ne \int_{1/2}^{1-\varepsilon }
\frac{t^{x-1}}{1-t}\ln^n t\,dt = \int_{1/2}^1
\frac{t^{x-1}}{1-t}\ln^n t\,dt$$ the integrals being convergent.
Further, from the Lemma 3.1 we see that the neutrix limit of the
function $$\sum_{i=0}^{r-1}(-1)^i \int_\varepsilon ^{1/2}
t^{x+i-1}\ln^n t\,dt $$ exists and implying that
$$\Ne \int_\varepsilon ^{1-\varepsilon
}\frac{t^{x-1}}{1-t}\ln^n t\,dt $$ exists. This proves the
existence of the function $\psi^{(n)}(x)$ for $n=0,1,2,\ldots,$
and all $x.$

\bigskip \noindent
Before giving our main theorem, we note that
$$\psi^{(n)}(x)=-\Ne \int_\varepsilon
^1\frac{t^{x-1}}{1-t}\ln^n t\,dt$$ since the integral is
convergent in the neighborhood of the point $t=1.$

\medskip \noindent
{\bf Theorem 3.4.} The function $\psi^{(n)}(-m)$ exists and \be
\psi^{(n)}(-m)= \sum_{i=1}^m\frac{n!}{i^{n+1}}+(-1)^{n+1}n!\zeta
(n+1)\ee for $n=1,2,\ldots$ and $m=0,1,2,\ldots,$ where $\zeta(n)$
denotes zeta function.

\medskip \noindent
{\bf Proof.} From Theorem 3.3, we have \beq \psi^{(n)}(-m)&=& -\Ne
\int_\varepsilon^1\frac{t^{-m-1}\ln ^n t}{1-t}\,dt \nonumber \\
&=&-\Ne \int_\varepsilon^1\Bigl [\sum_{i=1}^{m+1}
t^{-i}+(1-t)^{-1}\Bigr ]\ln ^nt\,dt. \eeq We first of all evaluate
the neutrix limit of integral $\int_\varepsilon^1 t^{-i}\ln^n
t\,dt$ for $i=1,2,\ldots$ and $n=1,2,\ldots.$

\noindent It follows from Lemma 3.2 that \be \Ne
\int_\varepsilon^1 t^{-i}\ln^n t\,dt=-\frac{n!}{(i-1)^{n+1}}\qquad
(i>1).\ee

\noindent For $i=1,$ we have
$$\int_\varepsilon^1 t^{-1}\ln^n t\,dt = O(\varepsilon ).$$
Next \beq \int_\varepsilon^1 \frac{\ln ^n t}{1-t}\,dt &=&\int_0^1
\frac{\ln ^n t}{1-t}\,dt=\sum_{k=0}^\infty \int_0^1 t^k\ln ^nt\,dt
\nonumber \\ &=&\sum_{k=0}^\infty \frac{(-1)^n n!}{(k+1)^{n+1}}
\nonumber \\ &=& (-1)^n n!\zeta(n+1), \eeq where $$\int_0^1 t^k\ln
^nt\,dt = \frac{(-1)^n n!}{(k+1)^{n+1}}.$$
It now follows from
equations (30), (31) and (32) that \beqa \psi^{(n)}(-m)&=& -\Ne
\int_\varepsilon^1\frac{t^{-m-1}\ln ^n t}{1-t}\,dt  \\
&=& -\sum_{i=1}^{m+1}\Ne \int_\varepsilon^1 t^{-i}\ln^n t\,dt-\Ne
\int_\varepsilon^1 \frac{\ln^n t}{(1-t)}\,dt. \\
&=& \sum_{i=1}^m\frac{n!}{i^{n+1}}+(-1)^{n+1}n!\zeta (n+1) \eeqa
implying equation (29).

\bigskip \noindent
Note that the digamma function $\psi(x)$ can be
defined by \be \psi (x)= -\gamma+\Ne \int_{\varepsilon}^1
\frac{1-t^{x-1}}{1-t}\,dt \ee for all $x.$

\bigskip \noindent Also note that using Lemma 3.1 we have \beqa
\int_{\varepsilon}^1 \frac{1-t^{-m-1}}{1-t}\,dt
&=&-\sum_{i=1}^{m+1}\int_{\varepsilon}^1
t^{-i}\,dt \\
&=&-[\ln 1-\ln
\varepsilon]-\sum_{i=2}^{m+1}\frac{[1-\varepsilon^{-i+1}]}{-i+1}
\eeqa and it follows from equation (33) that \beq \psi (-m) &=&
-\gamma+\Ne \int_{\varepsilon}^1
\frac{1-t^{-m-1}}{1-t}\,dt=-\gamma
+\sum_{i=1}^m i^{-1} \nonumber \\
&=&-\gamma +\phi (m) \eeq which was also obtained in
\cite{kn:fkurb2} and \cite{kn:biljana}. Equation (34) can be
regarded as the definition of the digamma function for negative
integers.

\medskip\noindent {\bf Emin Özçað } \\ Department of Mathematics \\
Hacettepe  University, Beytepe, Ankara, Turkey\\ e-mail:
ozcag1@hacettepe.edu.tr

\medskip \noindent
{\bf Ýnci Ege}\\ Department of Mathemacis \\ Adnan Menderes
University, Aydýn, Turkey \\
e-mail : iege@adu.edu.tr

\medskip\noindent


\bigskip
\begin{thebibliography}{99}
\bibitem{kn:Abr} M. Abramowitz and I. Stegun, Hand of Mathematical
Functions, Dover Publications, Newyork, 1965.
\bibitem{kn:c} J.G. van der Corput, Introduction to the neutrix
calculus, J. Analyse Math., {\bf 7}(1959), 291-398.
\bibitem{kn:f1} B. Fisher, On defining the incomplete gamma function $\gamma
(-m,x_-),$ Integral Trans. Spec. Funct., 15(6)(2004), 467-476.
\bibitem{kn:fbkil} B. Fisher, B. Jolevska-Tuneska and A. Kýlýçman, On
defining the incomplete gamma function, Integral Trans. Spec.
Funct., 14(4)(2003), 293-299.
\bibitem{kn:fkurb2} Fisher, B. and Kuribayashi, Y., Some results on the
Gamma function, J. Fac. Ed. Tottori Univ. Mat. Sci., 3(2)(1988)
111-117.
\bibitem{kn:gautschi} W. Gautschi, the incoplete gamma function
since Tricomi. Tricomi's Ideas Contemporary Applied Mathematics.
Number 147 in Atti dei Convegni Lincei. Roma : Accedemia nazionale
dei Lincel, Roma, Italy, 1998, 203-237.
\bibitem{kn:gs}I. M. Gel'fand and G. E. Shilov,
Generalized Functions, Vol. I., Academic Press, Newyork/London,
(1964).
\bibitem{kn:rr} I. S. Gradshhteyn, and I. M. Ryzhik, Tables of integrals, Series, and
Products, Academic Press, San Diego, 2000.
\bibitem{kn:biljana} Jolevska-Tuneska, B. and Jolevski, I., Some
results on the digamma function, Appl. Math. Inf. Sci., 7(2013)
167-170.
\bibitem{kn:kol3} Kölbig, K. S., On the integral
$\int_0^1 x^{\nu -1}(1-x)^{-\lambda}\ln ^mx\,dx,$ J. Comput. Appl.
Math., 18(1987) 369-394.
\bibitem{kn:ozinci} E. Özçað, Ý. Ege, H. Gürçay and
B. Jolevska-Tuneska, Some remarks on the incomplete gamma
function, in: Kenan Ta\c s et al. (Eds.), Mathematical Methods in
Engineering, Springer, Dordrecht, 2007, pp. 97-108.
\bibitem{kn:ozinci2} Özçað, E., Ege, Ý., Gürçay, H. and
Jolevska-Tuneska, B., On partial derivatives of the incomplete
beta function, Appl. Math. Lett., 21(2008) 675-681.
\bibitem{kn:Thomp1} I. Thompson, A note on real zeros of the
incomplete gamma function, Integral Trans. Spec. Funct.,
23(6)(2012), 445-453.
\bibitem{kn:Thomp2} I. Thompson, Algorithm 926: Incomplete Gamma
Function with Negative Arguments, ACM Trans. Math. Softw.
39(2)(2013), 9 pp.
\bibitem{kn:Win} S. Winitzki, computing the incomplete gamma
function to arbitrary precision, In Proceedings of International
Conference on Computational Science and Its Applications
(ICCSA'03). Lecturer Notes in Computer Science, Vol. 2667 (2003),
790-798.

\end{thebibliography}
\end{document}